\documentclass[a4paper,11pt,english]{smfart}

\usepackage{amssymb}
\usepackage{url}
\usepackage[T1]{fontenc}
\RequirePackage{calrsfs}
\DeclareSymbolFont{rsfscript}{OMS}{rsfs}{m}{b}
\DeclareSymbolFontAlphabet{\mathrsfs}{rsfscript}
\usepackage[utopia,expert]{mathdesign}
\usepackage{lscape}

\usepackage{enumitem}

\usepackage[leftbars]{changebar}

\usepackage{palatino}
\usepackage{rotating}
\usepackage{graphicx}

\usepackage{enumerate}

\usepackage{color}
\definecolor{shadecolor}{gray}{0.90}

\def\bfit{\bfseries\itshape}

\def\proofname{D\'emonstration}

\input xypic
\xyoption{all}
\xyoption{arc}

\makeindex

\def\indexname{Index des notations}
\newtheorem{theo}{Theorem}[section]
\def\thetheo{\thesection.\arabic{theo}}
\newtheorem{prop}[theo]{Proposition}

\newtheorem{coro}[theo]{Corollary}

\newtheorem{conj}[theo]{Conjecture}

\renewcommand\theequation{\thetheo}
\def\equat{\refstepcounter{theo}\begin{equation}}
\def\endequat{\end{equation}}

\renewcommand\thetable{\thetheo}

\renewcommand\thesection{\arabic{section}}
\renewcommand\thesubsection{\thesection.\Alph{subsection}}

\def\contentsname{Table des mati\`eres}
\def\listfigurename{Liste des figures}
\def\listtablename{Liste des tables}
\def\appendixname{Appendice}

\def\refname{R\'ef\'erences}

    \def\CM{{\mathbb{C}}}

  \def\hG{{\mathfrak h}}

  \def\lG{{\mathfrak l}}  
  \def\mG{{\mathfrak m}}  
    \def\NM{{\mathbb{N}}}
    
  \def\pG{{\mathfrak p}}

\def\SG{{\mathfrak S}}  \def\sG{{\mathfrak s}}  
  \def\tG{{\mathfrak t}}

    \def\ZM{{\mathbb{Z}}}

  \def\mGB{{\boldsymbol{\mathfrak m}}}

  \def\pGB{{\boldsymbol{\mathfrak p}}}

  \def\tGB{{\boldsymbol{\mathfrak t}}}

  \def\eb{{\mathbf e}}  \def\EC{{\mathcal{E}}}
    \def\FC{{\mathcal{F}}}
\def\Gb{{\mathbf G}}    
\def\Hb{{\mathbf H}}    \def\HC{{\mathcal{H}}}
  \def\ib{{\mathbf i}}

\def\Lb{{\mathbf L}}    \def\LC{{\mathcal{L}}}

    \def\OC{{\mathcal{O}}}
\def\Pb{{\mathbf P}}  \def\pb{{\mathbf p}}  \def\PC{{\mathcal{P}}}

\def\Sb{{\mathbf S}}    
\def\Tb{{\mathbf T}}

\def\Zb{{\mathbf Z}}

    \def\CCB{{\boldsymbol{\mathcal{C}}}}
    \def\DCB{{\boldsymbol{\mathcal{D}}}}
\def\Erm{{\mathrm{E}}}

\def\Hrm{{\mathrm{H}}}    \def\HCB{{\boldsymbol{\mathcal{H}}}}

\def\Krm{{\mathrm{K}}}

\def\Srm{{\mathrm{S}}}    \def\SCB{{\boldsymbol{\mathcal{S}}}}
\def\Trm{{\mathrm{T}}}

    \def\XCB{{\boldsymbol{\mathcal{X}}}}
    \def\YCB{{\boldsymbol{\mathcal{Y}}}}
\def\Zrm{{\mathrm{Z}}}    \def\ZCB{{\boldsymbol{\mathcal{Z}}}}

\def\a{\alpha}
\def\b{\beta}
\def\g{\gamma}

\def\e{\varepsilon}
\def\ph{\varphi}
\def\l{\lambda}
\def\L{\Lambda}

\def\o{\omega}
\def\O{\Omega}

\def\th{\theta}

\def\z{\zeta}

\def\Omeb{{\boldsymbol{\Omega}}}

           \def\phh{{\hat{\varphi}}}

\DeclareMathOperator{\Hom}{{\mathrm{Hom}}}

\DeclareMathOperator{\Id}{{\mathrm{Id}}}

\DeclareMathOperator{\im}{{\mathrm{Im}}}

\DeclareMathOperator{\Ind}{{\mathrm{Ind}}}

\DeclareMathOperator{\Irr}{{\mathrm{Irr}}}
\DeclareMathOperator{\Ker}{{\mathrm{Ker}}}

\DeclareMathOperator{\Res}{{\mathrm{Res}}}

\DeclareMathOperator{\Tr}{{\mathrm{Tr}}}

\def\to{\rightarrow}
\def\longto{\longrightarrow}
\def\injto{\hookrightarrow}

\def\fonction#1#2#3#4#5{\begin{array}{rccc}
{#1} : & {#2} & \longto & {#3}  \\
& {#4} & \longmapsto & {#5} 
\end{array}}

\def\fonctiol#1#2#3#4{\begin{array}{ccl}
{#1} & \longto & {#2} \\
{#3} & \longmapsto & {#4} 
\end{array}}

\def\DS{\displaystyle}

\def\finl{~$\blacksquare$}

\def\lexp#1#2{\kern\scriptspace\vphantom{#2}^{#1}\kern-\scriptspace#2}
\def\le{\hspace{0.1em}\mathop{\leqslant}\nolimits\hspace{0.1em}}
\def\ge{\hspace{0.1em}\mathop{\geqslant}\nolimits\hspace{0.1em}}

\mathchardef\inferieur="321E
\mathchardef\superieur="321F

\def\eqna{\begin{eqnarray*}}
\def\endeqna{\end{eqnarray*}}

\def\itemth#1{\item[${\mathrm{(#1)}}$]}

\catcode`\@=11
\long\def\@car#1#2\@nil{#1}
\long\def\@first#1#2{#1}
\long\def\@second#1#2{#2}
\long\def\ifempty#1{\expandafter\ifx\@car#1@\@nil @\@empty
  \expandafter\@first\else\expandafter\@second\fi}
\catcode`\@=12

\def\Hbov{{\bar{\Hb}}}
\def\Zbov{{\bar{\Zb}}}
\def\mGBov{{\bar{\mGB}}}

\def\GL{{\mathrm{GL}}}

\DeclareMathOperator{\Ref}{Ref}

\def\cow{{\mathrm{co}(W)}}

\theoremstyle{remark}
\newtheorem{rema}[theo]{Remark}

\newtheorem{exemple}[theo]{Example}

\theoremstyle{plain}

\def\BIL{LR}
\def\GAUCHE{L}
\def\CAR{CAR}
\def\FAM{FAM}

\def\grad{{\mathrm{gr}}}

\def\calo{{\Crm\Mrm}}

\def\xyinj{\ar@{^{(}->}}
\def\xysur{\ar@{->>}}

\def\isomorphisme#1{{\boldsymbol{[}}\hskip0.5mm #1\hskip0.5mm {\boldsymbol{]}}}

\bigskip

\makeatletter
\def\hlinewd#1{%
\noalign{\ifnum0=`}\fi\hrule \@height #1 %
\futurelet\reserved@a\@xhline}
\makeatother

\newlength\epaisLigne

\usepackage{multirow}

\newcommand{\longiso}{\stackrel{\sim}{\longrightarrow}}

\def\modules{\operatorname{\!-mod}\nolimits}

\def\gr{\operatorname{gr}\nolimits}
\def\Rees{\operatorname{Rees}\nolimits}
\def\codim{\operatorname{codim}\nolimits}

\makeatletter
\def\hlinewd#1{%
\noalign{\ifnum0=`}\fi\hrule \@height #1 %
\futurelet\reserved@a\@xhline}
\makeatother

\usepackage{multirow}

\newcommand{\longinjto}{\lhook\joinrel\longrightarrow}
\newcommand{\longsurto}{\relbar\joinrel\twoheadrightarrow}

\usepackage{lscape}

\usepackage{pstricks}

\def\GL{\operatorname{\Gb\Lb}\nolimits}

\def\pt{{\mathbf{pt}}}

\def\cow{{{\mathrm{co}}(W)}}
\def\gr{{\mathrm{gr}}}

\begin{document}

%\baselineskip=16pt
%\large\baselineskip=20pt
%\Large\baselineskip=24pt

\title{On the cohomology of Calogero-Moser spaces}

\author{{\sc C\'edric Bonnaf\'e}}
\address{
Institut Montpelli\'erain Alexander Grothendieck (CNRS: UMR 5149), 
Universit\'e Montpellier 2,
Case Courrier 051,
Place Eug\`ene Bataillon,
34095 MONTPELLIER Cedex,
FRANCE} 

\makeatletter
\email{cedric.bonnafe@umontpellier.fr}
\makeatother

\author{{\sc Peng Shan}}

\address{Yau Mathematical Sciences Center, Tsinghua University, 100086, BEIJING, CHINA}
\email{pengshan@tsinghua.edu.cn}

\date{\today}

\thanks{The first author is partly supported by the ANR 
(Project No ANR-16-CE40-0010-01 GeRepMod). \\
The second author is partly supported by the ANR (Project No ANR-13-8501-001-01 VARGEN).}

\maketitle
\pagestyle{myheadings}
\markboth{\sc C. Bonnaf\'e \& P. Shan}{\sc Cohomology of Calogero-Moser spaces}

\begin{abstract} 
We compute the equivariant cohomology of smooth Calogero-Moser spaces and 
some associated symplectic resolutions of symplectic quotient singularities.
\end{abstract}

\section{Notation and main results}\label{sec:notations}

\medskip

Throughout this note, we will abbreviate $\otimes_\CM$ as $\otimes$. By an algebraic variety, 
we mean a reduced scheme of finite type over $\CM$. 

\bigskip
\def\cod{{\mathrm{cod}}}

\subsection{Reflection group} 
Let $V$ be a $\CM$-vector space of finite dimension $n$ and let $W$ be a finite 
subgroup of $\GL_\CM(V)$. We set
$$\Ref(W)=\{s \in W~|~\dim_\CM V^s=n-1\}$$
and we assume that
$$W=\langle \Ref(W) \rangle.$$
We set $\e : W \to \CM^\times$, $w \mapsto \det(w)$. 

If $s \in \Ref(W)$, we denote by $\a_s^\vee$ and $\a_s$ two elements of $V$ and $V^*$, respectively, 
such that $V^s=\Ker(\a_s)$ and $V^{* s}=\Ker(\a_s^\vee)$, where $\a_s^\vee$ is viewed as a linear 
form on $V^*$.

If $w \in W$, we set
$$\cod(w)=\codim_\CM(V^w)$$
and we define a filtration $\FC_\bullet(\CM W)$ of the group algebra of $W$ as follows: let
$$\FC_i(\CM W)=\bigoplus_{\cod(w) \le i} \CM w.$$
Then 
$$\CM\Id_V=\FC_0(\CM W)\subset \FC_1(\CM W)\subset \cdots\subset\FC_n(\CM W)=\CM W=\FC_{n+1}(\CM W)= \cdots$$ 
is a filtration of $\CM W$. For any subalgebra $A$ of $\CM W$, we set $\FC_i(A)=A\cap\FC_i(\CM W)$, 
so that 
$$\CM\Id_V=\CM\FC_0(A)\subset \FC_1(A)\subset \cdots\subset\FC_n(A)=A=\FC_{n+1}(A)=\cdots$$ 
is also a filtration of $A$. Let $\hbar$ be a formal variable, and write 
\begin{align*}
&\Rees_\FC^\bullet(A)=\bigoplus_{i\geqslant 0}\hbar^i\FC_i(A)\subset \CM[\hbar] \otimes A
\qquad(\text{the {\it Rees algebra}}),\\
&\gr_\FC^\bullet(A)=\bigoplus_{i\geqslant 0}\FC_i(A)/\FC_{i-1}(A).
\end{align*}
Recall that $\gr_\FC^\bullet(A) \simeq \Rees_\FC^\bullet(A)/\hbar \Rees_\FC^\bullet(A)$. 

\bigskip
\def\la{\langle}
\def\ra{\rangle}

\subsection{Rational Cherednik algebra at $t=0$} 
Throughout this note, we fix a function $c : \Ref(W) \to \CM$ which 
is invariant under conjugacy. We define the  rational Cherednik algebra $\Hb_c$ to be the quotient 
of the algebra $\Trm(V\oplus V^*)\rtimes W$ (the semi-direct product of the tensor algebra 
$\Trm(V \oplus V^*)$ with the group $W$) 
by the relations 
$$\begin{cases}
[x,x']=[y,y']=0,\\
[y,x]=\DS{\sum_{s\in\Ref(W)}(\e(s)-1)c_s
\frac{\la y,\alpha_s\ra\la\alpha_s^\vee,x\ra}{\la\alpha_s^\vee,\alpha_s\ra}s},
\end{cases}\leqno{(\HC_c)}$$
for all $x$, $x'\in V^*$, $y$, $y'\in V$. Here $\la\ ,\ \ra: V\times V^*\to\CM$ is the standard pairing. 
The first commutation relations imply that 
we have morphisms of algebras $\CM[V] \to \Hb_c$ and $\CM[V^*] \to \Hb_c$. 
Recall~\cite[Theorem~1.3]{EG} 
that we have an isomorphism of $\CM$-vector spaces 
\equat\label{eq:pbw}
\CM[V] \otimes \CM W \otimes \CM[V^*] \longiso \Hb_c
\endequat
induced by multiplication (this is the so-called {\it PBW-decomposition}). 

We denote by $\Zb_c$ the center of $\Hb_c$: it is well-known~\cite{EG} that 
$\Zb_c$ is an integral domain, which is integrally closed and contains 
$\CM[V]^W$ and $\CM[V^*]^W$ as subalgebras (so it contains $\Pb=\CM[V]^W \otimes \CM[V^*]^W$), 
and that it is a free $\Pb$-module of rank $|W|$. We denote by $\ZCB_{\! c}$ the 
affine algebraic variety whose ring of regular functions $\CM[\ZCB_{\! c}]$ is $\Zb_c$: 
this is the {\it Calogero-Moser space} associated with the datum $(V,W,c)$. 

Using the PBW-decomposition, we define a $\CM$-linear map 
$\Omeb^c_\Hb : \Hb_c \longto \CM W$
by 
$$\Omeb^c_\Hb(f w g)=f(0)g(0)w$$
for all $f \in \CM[V]$, $g \in \CM[V^*]$ and $w \in \CM W$. This map is $W$-equivariant 
for the action on both sides by conjugation, so it induces a well-defined $\CM$-linear map 
$$\Omeb^c : \Zb_c \longto \Zrm(\CM W).$$
Recall from~\cite[Corollary~4.2.11]{BR} that $\Omeb^c$ is a morphism of algebras, and that 
\equat\label{eq:lissite-omega}
\text{\it $\ZCB_{\! c}$ is smooth if and only if $\Omeb^c$ is surjective.}
\endequat
The ``only if'' part is essentially due to Gordon~\cite[Corollary~5.8]{gordon} 
(see also~\cite[Proposition~9.6.6~and~(16.1.2)]{BR}) while 
the ``if'' part follows from the work of Bellamy, Schedler and Thiel~\cite[Corollary~1.4]{BST}.

\subsection{Grading}
The algebra $\Trm(V\oplus V^*)\rtimes W$ can be $\ZM$-graded in such a way that the generators have the following degrees
$$
\begin{cases}
\deg(y)=-1 & \text{if $y \in V$,}\\
\deg(x)=1 & \text{if $x \in V^*$,}\\
\deg(w)=0 & \text{if $w \in W$.}
\end{cases}
$$
This descends to a $\ZM$-grading on $\Hb_c$, because the defining relations $(\HC_c)$ are homogeneous. Since the center of a graded algebra is always graded, the subalgebra $\Zb_c$ is also $\ZM$-graded.  So the Calogero-Moser space $\ZCB_{\! c}$ 
inherits a regular $\CM^\times$-action. Note also that
by definition $\Pb=\CM[V]^W \otimes \CM[V^*]^W$ is clearly a graded subalgebra of $\Zb_c$.

\medskip

\subsection{Main results}\label{sub:main}
For a complex algebraic variety $\XCB$ (equipped with its classical topology), 
we denote by $\Hrm^i(\XCB)$ its $i$-th singular cohomology group with coefficients in $\CM$. If 
$\XCB$ carries a regular action of a torus $\Tb$, 
we denote by $\Hrm_\Tb^i(\XCB)$ its $i$-th $\Tb$-equivariant cohomology group 
(still with coefficients in $\CM$). Note that 
$\Hrm^{2\bullet}(\XCB)=\bigoplus_{i \ge 0} \Hrm^{2i}(\XCB)$ is a graded $\CM$-algebra and 
$\Hrm^{2\bullet}_\Tb(\XCB)=\bigoplus_{i \ge 0} \Hrm^{2i}_\Tb(\XCB)$ 
is a graded $\Hrm^{2\bullet}_\Tb(\pt)$-algebra. 
The following 
result~\cite[Theorem~1.8]{EG} describes the algebra structure on the cohomology of $\ZCB_{\! c}$
(with coefficients in $\CM$):

\medskip

\begin{theo}[Etingof-Ginzburg]\label{theo:eg}
Assume that $\ZCB_{\! c}$ is smooth. Then:
\begin{itemize}
\itemth{a} $\Hrm^{2i+1}(\ZCB_{\! c})=0$ for all $i$. 

\itemth{b} There is an isomorphism of graded $\CM$-algebras 
$\Hrm^{2\bullet}(\ZCB_{\! c}) \simeq \gr_\FC^\bullet (\Zrm(\CM W))$.
\end{itemize}
\end{theo}

\medskip

In this note, we prove an equivariant version of this statement (we identify 
$\Hrm^*_{\CM^\times}(\pt)$ with $\CM[\hbar]$ in the usual way):

\medskip

\noindent{\bfit Theorem A. ---} 
{\it Assume that $\ZCB_{\! c}$ is smooth. Then:
\begin{itemize}
\itemth{a} $\Hrm^{2i+1}_{\CM^\times}(\ZCB_{\! c})=0$ for all $i$.

\itemth{b} There is an isomorphism of graded $\CM[\hbar]$-algebras 
$\Hrm^{2\bullet}_{\CM^\times}(\ZCB_{\! c}) \simeq \Rees_\FC^\bullet(\Zrm(\CM W))$.
\end{itemize}}

\medskip

Note that Theorem~A(a) just follows from the statement~(a) of Etingof-Ginzburg Theorem 
by Proposition~\ref{prop:localisation}(a) below. 
As a partial consequence of Theorem~A, we also obtain the following application to the equivariant 
cohomology of some symplectic resolutions.

\bigskip
\noindent{\bfit Theorem B. ---} 
{\it Assume that the symplectic quotient singularity $(V \times V^*)/W$ admits 
a symplectic resolution $\pi : \XCB \longto (V \times V^*)/W$. Recall that the 
$\CM^\times$-action on $(V \times V^*)/W$ lifts (uniquely) to $\XCB$ 
(see~\cite[Theorem~1.3(ii)]{K}). Then:
\begin{itemize}
\itemth{a} $\Hrm^{2i+1}_{\CM^\times}(\XCB)=0$ for all $i$.

\itemth{b} There is an isomorphism of graded $\CM[\hbar]$-algebras 
$\Hrm^{2\bullet}_{\CM^\times}(\XCB) \simeq \Rees_\FC^\bullet(\Zrm(\CM W))$.
\end{itemize}}

\medskip

Note that for $W=\SG_n$ acting on $\CM^n$, Theorem~B describes the equivariant 
cohomology of the Hilbert scheme of $n$ points in $\CM^2$: this was already 
proved by Vasserot~\cite{vasserot}. 
In~\cite[Conjecture~1.3]{GK}, Ginzburg-Kaledin proposed 
a conjecture for the equivariant cohomology of a symplectic resolution 
of a symplectic quotient singularity $E/G$, where $E$ is a finite dimensional 
symplectic vector space and $G$ is a finite subgroup of $\Sb\pb(E)$. 
However, their conjecture cannot hold as stated, because they considered the $\CM^\times$-action by dilatation, which is contractible. 
Theorem~B shows that the correct equivariant cohomological realization of the Rees algebra is provided by the symplectic $\CM^\times$-action, 
which exists only when $G$ stabilizes a Lagrangian subspace of $E$. 

\bigskip

\begin{rema}\label{rem:classification} 
Recall from the works of Etingof-Ginzburg~\cite{EG}, 
Ginzburg-Kaledin~\cite{GK}, Gordon~\cite{gordon} and Bellamy~\cite{bellamy} that 
the existence of a symplectic resolution of $(V \times V^*)/W$ is equivalent to the existence 
of a parameter $c$ such that $\ZCB_c$ is smooth, and that it can only occur 
if all the irreducible components of $W$ are of type $G(d,1,n)$ (for some $d$, $n \ge 1$) 
or $G_4$ in Shephard-Todd classification.\finl 
\end{rema}

\bigskip

\definecolor{shadecolor}{gray}{0.95}
\definecolor{boite}{gray}{0.98}
\def\springer#1{\begin{centerline}{\fcolorbox{black}{shadecolor}{~
    \begin{minipage}[t]{0.8\textwidth}{\vphantom{~}{\itshape #1}\vphantom{$A_{\DS{A_A}}$}}
            \end{minipage}~}}\end{centerline}\medskip}

\subsection{Conjectures}\label{sub:conjectures}

In~\cite[Chapter~16,~Conjectures~COH~and~ECOH]{BR}, Rouquier and the first author proposed the 
following conjecture which aims to generalize the above Etingof-Ginzburg Theorem~\ref{theo:eg} into two 
directions: it includes singular Calogero-Moser spaces and it extends to equivariant cohomology. 

\medskip

\begin{conj}\label{conj:h}
With the above notation, we have:
\begin{itemize}%[leftmargin=1cm]
\itemth{1} $\Hrm^{2i+1}(\ZCB_{\! c})=0$ for all $i$.

\itemth{2} We have an isomorphism of graded $\CM$-algebras 
$\Hrm^{2\bullet}(\ZCB_{\! c}) \simeq \gr_\FC(\im(\Omeb^c))$.

\itemth{3} We have an isomorphism of graded $\CM[\hbar]$-algebras 
$\Hrm^{2\bullet}_{\CM^\times}(\ZCB_{\! c}) \simeq \Rees_\FC(\im(\Omeb^c))$.
\end{itemize}
\end{conj}

\medskip

% 
% This conjecture is then 
% extended to the equivariant setting~\cite[Chapter~17,~Conjecture~ECOH]{BR}. 
% %and 
% %Theorem~A above is a particular case of this last conjecture 
% %(see Section~\ref{sub:conjectures} for precise statements). 
% 
% %Recall from~\cite[Chapter~17]{BR} the following conjectures.
% 
% %We denote by $\hbar$ the first Chern class of the $\CM^\times$-module of weight $1$ in 
% %$\Hrm^2_{\CM^\times}(\pt)$, so that $\hbar$ is an indeterminate over $\CM$ and
% %Recall that $\Hrm^{\bullet}_{\CM^\times}(\pt)=\CM[\hbar]$, with $\hbar$ an indeterminate of degree $2$. 
% %Denote this ring by $R$.
% \medskip
% 
% \begin{conj}[Equivariant cohomology]\label{conj:eh}
% With the above notation, we have:
% \begin{itemize}%[leftmargin=1.3cm]
% \itemth{1} $\Hrm^{2i+1}_{\CM^\times}(\ZCB_{\! c})=0$ for all $i$. 
% 
% 
% \end{itemize}
% \end{conj}
% 
% \medskip

By~(\ref{eq:lissite-omega}), when $\ZCB_{\! c}$ is smooth, 
the image of $\Omeb^c$ coincide with the center of $\CM W$. 
So Theorem A proves this conjecture for smooth $\ZCB_{\! c}$.

\medskip

%\begin{theo}[Etingof-Ginzburg]
%If $\ZCB_{\! c}$ is smooth, then Conjecture~\ref{conj:h} holds.
%\end{theo}

%Our aim in this note is to prove the following equivariant version of this theorem
%\bigskip

%\begin{theo}\label{theo:equiv}
%If $\ZCB_{\! c}$ is smooth, then Conjecture~\ref{conj:eh} holds.
%\end{theo}

%\bigskip

We will also prove in Example~\ref{ex:cyclique} the following result:

\bigskip

\begin{prop}\label{prop:cyclique}
If $\dim_\CM(V)=1$, then Conjecture~\ref{conj:h} holds.
\end{prop}

\subsection{Structure of the paper}
The paper is organized as follows. The proof of Theorem~A relies on classical theorems on 
restriction to fixed points in cohomology and K-theory. 
In Section~\ref{sec:localisation}, we first recall basic properties on equivariant cohomology and 
equivariant K-theory and restriction to fixed points. Section~\ref{sec:local-cm} 
explains how these general principles can be applied to Calogero-Moser spaces. 
Theorem~A will be proved in Section~\ref{sec:proof-A}. The cyclic group case (Proposition~\ref{prop:cyclique}) 
will be handled in Example~\ref{ex:cyclique}. 
The proof of Theorem~B will be given in Section~\ref{sec:proof-B}.

\bigskip

\section{Equivariant cohomology, $\Krm$-theory and fixed points}\label{sec:localisation}

\medskip
\def\chern{\operatorname{ch}\nolimits}

\subsection{Equivariant cohomology}\label{sub:H}
Let $\XCB$ be a complex algebraic variety equipped with a regular action of a torus $\Tb$. 
Recall that the \emph{equivariant cohomology} of $\XCB$ is defined by
$$\Hrm^\bullet_\Tb(\XCB)=\Hrm^\bullet(E_{\Tb}\times_{B_\Tb}\XCB),$$
where $E_\Tb\to B_\Tb$ is a universal $\Tb$-bundle.
The pullback of the structural morphism $x: \XCB\to\pt$ yields a ring homomorphism 
$\Hrm^\bullet_\Tb(\pt)\to \Hrm^\bullet_\Tb(\XCB)$, which makes $\Hrm^\bullet_\Tb(\XCB)$ a graded $\Hrm^\bullet_\Tb(\pt)$-algebra.

Denote by $X(\Tb)$ the character lattice of $\Tb$. Then for each $\chi\in X(\Tb)$, denote by 
$\CM_\chi$ the one dimensional $\Tb$-module of character $\chi$, the first Chern class $c_\chi$ 
of the line bundle $E_{\Tb}\times_{\Tb}\CM_{\chi}$ on $B_{\Tb}$ is an element in $\Hrm^2(B_\Tb)$. 
Identify the vector space $\CM \otimes_\ZM X(\Tb)$ with the dual $\tGB^\ast$  of the Lie algebra $\tGB$ of 
$\Tb$ via $\chi\mapsto d\chi$. Then the assignment $\chi\mapsto c_\chi$ yields an isomorphism 
of graded $\CM$-algebras $S(\tGB^\ast)=\Hrm^{2\bullet}_\Tb(\pt)$.

\subsection{Equivariant $\Krm$-theory}\label{sub:K}
We denote by $\Krm_\Tb(\XCB)$ the Grothendieck ring of the category of 
$\Tb$-equivariant vector bundles on $\XCB$. Note that a $\Tb$-equivariant vector bundle on a point is the same as a finite dimensional $\Tb$-module. 
We have a canonical isomorphism $\Krm_\Tb(\pt)=\ZM[X(\Tb)]$ which sends the class of a $\Tb$-module $M$ to 
\equat\label{eq:dim-graduee}
\dim^\Tb(M)=\sum_{\chi\in X(\Tb)}\dim_\CM(M_\chi)\chi,
\endequat
where $M_\chi$ is the $\chi$-weight space in $M$.

Let $\hat{\Hrm}^{2\bullet}_\Tb(\XCB)$ 
be the completion of $\Hrm_\Tb^{2\bullet}(\XCB)$ 
with respect to the ideal $\bigoplus_{i > 0} \Hrm^{2i}_\Tb(\XCB)$. 
The {\it equivariant Chern character} provides a ring homomorphism
$$\chern_\XCB : \Krm_\Tb(\XCB) \longto \hat{\Hrm}_\Tb^{2\bullet}(\XCB)$$ 
with the following properties. 
First, when $\XCB$ is a point $\pt$, we have
\equat\label{eq:chernpt}
\fonction{\chern_\pt}{\Krm_\Tb(\pt)=\ZM[X(\Tb)]}{\hat{\Hrm}_\Tb^{2\bullet}(\pt)=\hat{\Srm}(\tGB^*)}{\chi}{\exp(d\chi).}
\endequat
Next the Chern character commutes with pullback. More precisely,
if $\YCB$ is another variety with a regular action of the same torus $\Tb$ 
and if $\ph : \XCB \to \YCB$ is a $\Tb$-equivariant morphism, then 
%pull-back by 
%$\ph$ induces two morphisms of algebras 
%$$\ph^*_\Krm : \Krm_\Tb(\YCB) \longto \Krm_\Tb(\XCB)$$
%$$\ph^*_\Hrm : \Hrm_\Tb^{2\bullet}(\YCB) \longto \Hrm_\Tb^{2\bullet}(\XCB)\leqno{\text{and}}$$
the following diagram commutes
\equat\label{eq:chern}
\diagram
\Krm_\Tb(\YCB) \rrto^{\DS{\ph^*}} \ddto_{\DS{\chern_\YCB}}&& \Krm_\Tb(\XCB) \ddto^{\DS{\chern_\XCB}}\\
&&\\
\hat{\Hrm}_\Tb^{2\bullet}(\YCB) \rrto^{\DS{\phh^*}} &&\hat{\Hrm}_\Tb^{2\bullet}(\XCB).
\enddiagram
\endequat
Here $\ph^*$ denotes both the pullback map in $\Krm$-theory or in equivariant cohomology, and 
$\phh^*$ is the map induced after completion by the pullback map in equivariant cohomology. 
In particular, by applying the above diagram to $\YCB=\pt$ and the structural morphism $\XCB \to \pt$,
we may view $\chern_\XCB$ as a morphism of algebras over $\Krm_\Tb(\pt)$, with the $\Krm_\Tb(\pt)$-algebra structure on $\hat{\Hrm}_\Tb^{2\bullet}(\XCB)$ provided by the embedding \eqref{eq:chernpt}.

%For instance, we may view $\chern_\XCB$ as a morphism of algebras over $\Krm_\Tb(\pt)$ 
%via the commutative diagram (here, $x : \XCB \to \pt$ is the structural morphism)
%$$\diagram
%\Krm_\Tb(\pt) \rrto^{\DS{x_\Krm^*}} \ddto_{\DS{\chern_\pt}}&& \Krm_\Tb(\XCB) \ddto^{\DS{\chern_\XCB}}\\
%&&\\
%\hat{\Hrm}_\Tb^{2\bullet}(\pt) \rrto^{\DS{\hat{x}_\Hrm^*}} &&\hat{\Hrm}_\Tb^{2\bullet}(\XCB)
%\enddiagram$$

\medskip
\def\local{i}

\subsection{Fixed points}
We denote by $\XCB^\Tb$ the (reduced) variety consisting of fixed points of $\Tb$ in $\XCB$.
Let $\local_\XCB : \XCB^\Tb \longinjto \XCB$ be the natural closed immersion. 
Since $\Tb$ acts trivially on $\XCB^\Tb$, we have 
$\Hrm^{\bullet}_\Tb(\XCB^\Tb) = \Hrm^{\bullet}_\Tb(\pt) \otimes \Hrm^{\bullet}(\XCB^\Tb)$ as $\Hrm^{\bullet}_\Tb(\pt)$-algebras. 
Recall that the $\Tb$-action on $\XCB$ is called \emph{equivariantly formal} if the Leray-Serre spectral sequence
$$\Erm^{pq}_2=\Hrm^p(B_\Tb;\Hrm^q(\XCB))\Longrightarrow \Hrm^{p+q}_{\Tb}(\XCB)$$
for the fibration $E_{\Tb}\times_\Tb \XCB\to B_\Tb$ degenerates at $E_2$.
We have the following standard result on equivariant cohomology (see for instance~\cite[Proposition~2.1]{GM}):
%the classical {\it localization theorem} (reference???):

\bigskip

\begin{prop}\label{prop:localisation}
Assume that $\Hrm^{2i+1}(\XCB)=0$ for all $i$. Then $\XCB$ is equivariantly formal, and:
\begin{itemize}
\itemth{a} There is an isomorphism of graded $\Hrm^\bullet_\Tb(\pt)$-modules 
$\Hrm^\bullet_\Tb(\XCB)\simeq \Hrm^\bullet_\Tb(\pt)\otimes \Hrm^\bullet(\XCB)$. 
In particular $\Hrm_\Tb^{2i+1}(\XCB)=0$ for all $i$.

\itemth{b} The pullback map
$\local_{\XCB}^* : \Hrm^{\bullet}_\Tb(\XCB) \to \Hrm_\Tb^{\bullet}(\XCB^\Tb)$ 
is injective.

\itemth{c} Let $\mG$ be the unique graded maximal ideal of $\Hrm_\Tb^\bullet(\pt)$. 
Then we have an isomorphism of algebras 
$\Hrm^{\bullet}(\XCB) \simeq \Hrm_{\Tb}^\bullet(\XCB) / \mG \Hrm_{\Tb}^\bullet(\XCB)$. 
\end{itemize}
\end{prop}

\bigskip

In particular, this shows that Conjectures~\ref{conj:h}(1) and~(3) imply 
Conjecture~\ref{conj:h}(2).

\bigskip

\begin{exemple}[Blowing-up]\label{ex:blowing}
Let $\YCB$ be an affine variety with a $\Tb$-action and let $\CCB$ be a 
$\Tb$-stable closed subvariety (not necessarily reduced) of $\YCB$. Let
$\XCB$ be the blowing-up of $\YCB$ along $\CCB$. Write $\pi : \XCB \to \YCB$ for
the natural morphism, and we equip $\XCB$ with the unique $\Tb$-action 
such that $\pi$ is equivariant. {\bfit We assume that $\XCB^\Tb$ is finite}. 
We write 
$$\Hrm^{2\bullet}_\Tb(\XCB^\Tb)=\bigoplus_{x \in \XCB^\Tb} \Srm(\tG^*) \eb_x,$$
where $\eb_x \in \Hrm_\Tb^0(\XCB^\Tb)$ is the primitive idempotent associated 
with $x$ (i.e., the fundamental class of $x$). 

Then $\DCB=\pi^*(\CCB)$ is a $\Tb$-stable effective Cartier divisor, and we denote by 
$\isomorphisme{\DCB}$ the class in $\Krm_\Tb(\XCB)$ of its associated line bundle 
(which is $\Tb$-equivariant). We denote by 
$\chern_\XCB^1(\isomorphisme{\DCB}) \in \Hrm_\Tb^2(\XCB)$ its first $\Tb$-equivariant 
Chern class. We want to compute $\local_{\XCB}^*(\chern^1_\Tb(\DCB))$. 

First, let $I$ be the ideal of $\CM[\YCB]$ associated with $\CCB$. As it is $\Tb$-stable, 
we can find a family of $\Tb$-homogeneous generators $(a_i)_{1 \le i \le k}$ of $I$. 
We denote by $\l_i \in X(\Tb)$ the $\Tb$-weight of $a_i$. The choice of 
this family of generators induces a $\Tb$-equivariant closed immersion  
$\XCB \injto \YCB \times \Pb^{k-1}(\CM)$. We denote by $\XCB_{\! i}$ the affine chart 
corresponding to ``$a_i \neq 0$''. If $x \in \XCB^\Tb$, we denote by 
$\ib(x) \in \{1,2,\dots,k\}$ an element such that $x \in \XCB_{\! \ib(x)}$. Then
\equat\label{eq:local-blowing}
\local_{\XCB}^*(\chern^1_\Tb(\DCB)) = -\hbar \sum_{x \in \DCB^\Tb} (d\l_{\ib(x)}) \eb_x.
\endequat
Indeed, we just need to compute the local equation of $\DCB$ around $x \in \XCB^\Tb$, 
and this can be done in $\XCB_{\! \ib(x)}$. But 
$\DCB \cap \XCB_{\! i}$ is principal for all $i$, defined by 
$$\DCB \cap \XCB_{\! i} = \{(y,\xi) \in \YCB \times \Pb^{k-1}(\CM)~|~(y,\xi) \in \XCB_{\! i}~\text{  and  }
a_i(y)=0\}.$$
So $\DCB \cap \XCB_{\! i}$ is defined by a $\Tb$-homogeneous equation of degree $\l_i \in X(\Tb)$, and 
so~(\ref{eq:local-blowing}) follows.\finl
\end{exemple}

\bigskip
\def\calo{{\mathrm{CM}}}

\section{Localization and Calogero-Moser spaces}\label{sec:local-cm}

\medskip

In this section, we apply the previous discussions to $\XCB=\ZCB_{\! c}$ and $\Tb=\CM^\times$. 
Denote by $q:\CM^\times\to \CM^\times$ the identity map. Then $X(\CM^\times)=q^\ZM$, and we have
$$\Krm_{\CM^\times}(\pt) = \ZM[q,q^{-1}],\qquad \Hrm_{\CM^\times}^{2\bullet}(\pt)=\CM[\hbar],$$ with $\hbar=c_q$, 
following the notation of Section \ref{sub:H}.
So $\hat{\Hrm}_{\CM^\times}^{2\bullet}(\pt)=\CM[[\hbar]]$ and the Chern map in this case is given by
$$\chern_\pt : \ZM[q,q^{-1}] \injto \CM[[\hbar]],\qquad q\mapsto \exp(\hbar).$$
Note also that a finite $\CM^\times$-module is nothing but a finite dimensional $\ZM$-graded vector 
space $M=\bigoplus_{i \in \ZM} M_i$ such that $\CM^\times$ acts on $M_i$ by the character $q^i$.
The identification $\Krm_{\CM^\times}(\pt) = \ZM[q,q^{-1}]$ sends the class of $M$ to its 
{\it graded dimension} (or {\it Hilbert series}):
$$\dim^\gr(M)=\sum_{i \in \ZM} \dim_\CM(M_i) q^i \in \NM[q,q^{-1}].$$

\subsection{Fixed points} 
For $\chi \in \Irr(W)$, we denote by $\o_\chi : \Zrm(\CM W) \to \CM$ the 
associated morphism of algebras, that is, $\o_\chi(z)=\chi(z)/\chi(1)$ is the scalar through which 
$z$ acts on the irreducible $\CM W$-module with the character $\chi$. We denote by $e_\chi^W$ 
the unique primitive idempotent of $\Zrm(\CM W)$ such that $\o_\chi(e_\chi^W)=1$. 
If $\EC$ is a subset of $\Irr(W)$, then we set $e_\EC^W=\sum_{\chi \in \EC} e_\chi^W$.

Now, consider the algebra homomorphism
$$\O^c_\chi = \o_\chi \circ \Omeb^c : \Zb_c \longsurto \CM.$$
Its kernel is a maximal ideal of $\Zb_c$: we denote by $z_\chi$ the corresponding closed point in $\ZCB_{\! c}$. 
It follows from~\cite[Lemma~10.2.3~and~(14.2.2)]{BR} that $z_\chi \in \ZCB_{\! c}^{\CM^\times}$ and that 
the map
\equat\label{eq:fixes}
\fonction{z}{\Irr(W)}{\ZCB^{\CM^\times}}{\chi}{z_\chi}
\endequat
is surjective. The fibers of this map are called the {\it Calogero-Moser $c$-families}.
They were first consider by Gordon~\cite{gordon} and Gordon-Martino~\cite{gordon-martino}: 
see also for instance~\cite[\S{9.2}]{BR}. Let $\calo_c(W)$ be the set of Calogero-Moser 
$c$-families. For $\EC\in \calo_c(W)$, we denote by $z_\EC\in \ZCB^{\CM^\times}$ its image under the map $z$. 
On the other hand, by~\cite[(16.1.2)]{BR} we have
\equat\label{eq:image-omega}
\im(\Omeb^c)=\bigoplus_{\EC \in \calo_c(W)} \CM e_\EC^W.
\endequat
Hence we get an isomorphism of $\CM$-algebras
\begin{align*}
\Hrm^{2\bullet}(\ZCB_{\! c}^{\CM^\times})\simeq \im(\Omeb^c),\qquad e_{z_\EC}\mapsto e_\EC^W,
\end{align*}
which extends to an isomorphism of $\CM[\hbar]$-algebras
$\Hrm^{2\bullet}_{\CM^\times}(\ZCB_{\! c}^{\CM^\times})\simeq\CM[\hbar] \otimes \im(\Omeb^c).$ 
For simplification, we set $\local_c = \local_{\ZCB_{\! c}} : \ZCB_{\! c}^{\CM^\times} \injto \ZCB$ and, 
under the above identification, we view the pullback map $\local_c^*$ 
as a morphism of algebras
$$\local_c^* : \Hrm^{2\bullet}_{\CM^\times}(\ZCB_{\! c}) \longto \CM[\hbar] \otimes \im(\Omeb^c).$$ 
So, by Proposition~\ref{prop:localisation}, Conjecture~\ref{conj:h} 
is implied by the following one:

\bigskip

\begin{conj}\label{conj:h-eh}
With the above notation, we have:
\begin{itemize}
\itemth{1} $\Hrm^{2i+1}(\ZCB_{\! c})=0$ for all $i$.

\itemth{2} $\im(\local_c^*) = \Rees_\FC(\im(\Omeb^c))$.
\end{itemize}
\end{conj}

\bigskip

\begin{rema}\label{rem:filtration}
% We assume in this remark that Conjecture~\ref{conj:h}(1) holds (i.e. that 
% $\Hrm^{2i+1}(\ZCB_{\! c})=0$ for all $i$). Then the Euler characteristic $\chi(\ZCB_{\! c})$ 
% of $\ZCB_{\! c}$ satisfies
% $$\chi(\ZCB_{\! c})=\sum_{i \ge 0} \dim_\CM \Hrm^{2i}(\ZCB_{\! c})$$
% $$\chi(\ZCB_{\! c})=\chi(\ZCB_{\! c}^{\CM^\times})=|\ZCB_{\! c}^{\CM^\times}|, \leqno{\text{and}}$$
% so Conjecture~\ref{conj:h} is compatible with these two equalities. 
%On the other hand, 
Set 
$$\FC_i^\Hrm(\im(\Omeb^c))=\{z \in \im(\Omeb^c)~|~\hbar^i z \in \im(\local_c^*)\}.$$
Then, by construction and Proposition~\ref{prop:localisation}, $\FC_\bullet^\Hrm(\im(\Omeb^c))$ 
is the filtration of $\im(\Omeb^c)$ satisfying
$$\Hrm^{2\bullet}_{\CM^\times}(\ZCB_{\! c}) \simeq \im(\local_c^*)=\Rees_{\FC^\Hrm}(\im(\Omeb^c))$$
$$\Hrm^{2\bullet}(\ZCB_{\! c}) \simeq \gr_{\FC^\Hrm}(\im(\Omeb^c)).\leqno{\text{and}}$$
So showing Conjecture~\ref{conj:h-eh} is then equivalent to showing that the filtrations 
$\FC_\bullet(\im(\Omeb^c))$ and $\FC_\bullet^\Hrm(\im(\Omeb^c))$ coincide. 

Note also that, since $\ZCB_{\! c}$ is an affine variety of dimension $2n$, we have 
$\Hrm^i(\ZCB_{\! c})=0$ for $i > 2n$, so this shows that 
$$\FC_n(\im(\Omeb^c)) = \im(\Omeb^c) = \FC_n^\Hrm(\im(\Omeb^c)).\leqno{(\flat)}$$
On the other hand, since $\ZCB_{\! c}$ is connected, we have
$$\FC_0(\im(\Omeb^c))=\FC_0^\Hrm(\im(\Omeb^c))=\CM.\leqno{(\sharp)}$$
These two particular cases will be used below.\finl 
\end{rema}

\bigskip

The following result, based on the previous remark, 
can be viewed as a reduction step 
for the proof of Conjecture~\ref{conj:h-eh}:

\bigskip

\begin{prop}\label{prop:reduction}
Assume that $\Hrm^{2i+1}(\ZCB_{\! c})=0$ and 
$$\dim_\CM(\Hrm^{2i}(\ZCB_{\! c}))=\dim_\CM(\FC_i(\im(\Omeb^c))/\FC_{i-1}(\im(\Omeb^c)))$$ 
for all $i$. Then 
Conjecture~\ref{conj:h-eh} holds if and only if 
$\Rees_\FC(\im(\Omeb^c)) \subset \im(\local_c^*)$. 
\end{prop}

\bigskip

\begin{proof}
Assume that $\Hrm^{2i+1}(\ZCB_{\! c})=0$ and 
$$\dim_\CM(\Hrm^{2i}(\ZCB_{\! c}))=\dim_\CM(\FC_i(\im(\Omeb^c))/\FC_{i-1}(\im(\Omeb^c)))$$ 
for all $i$. We keep the notation of Remark~\ref{rem:filtration}. 
It then follows from this remark and the hypothesis that 
$$\dim_\CM(\FC_i(\im(\Omeb^c))/\FC_{i-1}(\im(\Omeb^c)))=
\dim_\CM(\Hrm^{2i}(\ZCB_{\! c}))=\dim(\FC_i^\Hrm(\im(\Omeb^c))/\FC_{i-1}^\Hrm(\im(\Omeb^c)))$$
for all $i$. So, by induction, we get that $\dim_\CM(\FC_i(\im(\Omeb^c)))=\dim(\FC_i^\Hrm(\im(\Omeb^c)))$ 
for all $i$ (by the equality~($\sharp)$ of Remark~\ref{rem:filtration}). 
This shows that $\Rees_\FC(\im(\Omeb^c)) = \im(\local_c^*)$ if and only if 
$\Rees_\FC(\im(\Omeb^c)) \subset \im(\local_c^*)$, as desired. 
\end{proof}

\bigskip

\begin{exemple}\label{ex:cyclique}
Assume in this example, and only in this example, that $\dim_\CM(V)=1$. It is proved in \cite[Theorem~18.5.8]{BR} 
that in this case $\Hrm^{2i+1}(\ZCB_{\! c})=0$ and 
$$\dim_\CM(\Hrm^{2i}(\ZCB_{\! c}))=\dim_\CM(\FC_i(\im(\Omeb^c))/\FC_{i-1}(\im(\Omeb^c)))$$ 
for all $i$. %We keep the notation of the Remark~\ref{rem:filtration} and the proof 
%of Proposition~\ref{prop:reduction}. 
Since $\ZCB_{\! c}$ is affine of dimension $2$, 
we have $\Hrm^i(\ZCB_{\! c})=0$ if $i \not\in \{0,2\}$. So it follows from the equalities~($\flat$) 
and~($\sharp$) of Remark~\ref{rem:filtration} that Conjecture~\ref{conj:h-eh} holds in this case 
(this proves Proposition~\ref{prop:cyclique}).\finl
\end{exemple}

\bigskip

\subsection{Chern map} 
If $\EC$ is a Calogero-Moser family, we denote by $\mG_\EC \subset \Zb_c$ the ideal of functions vanishing 
at $z_\EC \in \ZCB_{\! c}^{\CM^\times}$. We also set 
$$\im(\Omeb^c)_\ZM=\bigoplus_{\EC \in \calo_c(W)} \ZM e_\EC^W.$$
We make the natural identification 
$\Krm_{\CM^\times}(\ZCB_{\! c}^{\CM^\times}) = \ZM[q,q^{-1}] \otimes_\ZM \im(\Omeb^c)_\ZM$. 
% and we denote by $\local_c^\Krm$ the map 
% $\local_{\ZCB_{\! c}}^* : \Krm_{\CM^\times}(\ZCB_{\! c}) \longto \Krm_{\CM^\times}(\ZCB_{\! c}^{\CM^\times})$. 
Through these identifications, the Chern map $\chern_{\ZCB_{\! c}^{\CM^\times}}$ just 
becomes the natural inclusion $\ZM[q,q^{-1}] \otimes_\ZM \im(\Omeb^c)_\ZM \injto \CM[[\hbar]] \otimes \im(\Omeb^c)$. 
Moreover, if $P$ is a $\ZM$-graded finitely generated projective $\Zb_c$-module, then the commutativity of the 
diagram~(\ref{eq:chern}) just says that
\equat\label{eq:local-chern}
\local_c^*(\chern_{\ZCB_{\! c}}(\isomorphisme{P})) = \sum_{\EC \in \calo_c(W)} 
\dim^\gr(P/\mG_\EC P) e_\EC^W \subset \CM[[\hbar]] \otimes \im(\Omeb^c).
\endequat

\bigskip

\section{Proof of Theorem~A}\label{sec:proof-A}

\medskip

\springer{\noindent{\bf Hypothesis and notation.} 
{\it We assume in this section, and only in this section, that 
$\ZCB_{\! c}$ is smooth.}}

\medskip

If $M=\bigoplus_{i \in \ZM} M_i$ and $N=\bigoplus_{i \in \ZM} N_i$ are two finite-dimensional 
graded $\CM W$-modules, we set
$$\langle M,N \rangle^\gr_W = \sum_{i,j \in \ZM} \la M_i,N_j \ra_W ~q^{i+j},$$
where $\la E , F \ra_W =\dim \Hom_{\CM W}(E,F)$ for any finite dimensional $\CM W$-modules $E$ and $F$. 
We extend this notation to the case where $M$ or $N$ is a graded virtual character 
(i.e. an element of $\ZM[q,q^{-1}] \otimes_\ZM \ZM\Irr(W)$).
Finally, we denote by $\CM[V]^\cow$ the {\it coinvariant algebra}, 
that is, the quotient of the algebra $\CM[V]$ by the 
ideal generated by the elements $f \in \CM[V]^W$ such that $f(0)=0$. Then 
$\CM[V]^\cow$ is a graded $\CM W$-module, which is isomorphic to the regular 
representation $\CM W$ when one forgets the grading. 

\medskip

\subsection{Localization in $\Krm$-theory}
Recall from~(\ref{eq:lissite-omega}) that 
the smoothness of $\ZCB_{\! c}$ implies that $\im(\Omeb^c)=\Zrm(\CM W)$, 
so that $z : \Irr(W) \to \ZCB_{\! c}^{\CM^\times}$ is bijective. 
Let $e=e_1^W=(1/|W|)\sum_{w \in W} w$. 
The smoothness of $\ZCB_{\! c}$ also implies that the functor
$$\fonctiol{\Hb_c\modules}{\Zb_c\modules}{M}{eM = e\Hb_c \otimes_{\Hb_c} M}$$
is an equivalence of categories~\cite[Theorem~1.7]{EG}. If $E$ is a finite dimensional 
{\it $\ZM$-graded} $\CM W$-module, the $\Hb_c$-module $\Hb_c \otimes_{\CM W} E$ 
is finitely generated, $\ZM$-graded and projective. Therefore, $e\Hb_c \otimes_{\CM W} E$ is a finitely generated 
graded projective $\Zb_c$-module, 
which can be viewed as a $\CM^\times$-equivariant vector bundle on $\ZCB_{\! c}$. 
For simplification, we set 
$$\chern_c(E)=\local_c^*(\chern_{\ZCB_{\! c}}(e\Hb_c \otimes_{\CM W} E)) 
\in \CM[[\hbar]] \otimes \im(\Omeb^c).$$ 

\bigskip

\begin{prop}\label{prop:chern}
Assume that $\ZCB_{\! c}$ is smooth. 
Let $E$ be a finite dimensional graded $\CM W$-module. Then
$$\chern_c(E) =\sum_{\chi \in \Irr(W)} 
\frac{\langle \chi , \CM[V]^\cow \otimes E \rangle_W^\gr}{\la \chi, \CM[V]^\cow \ra_W^\gr} e_\chi^W.$$
\end{prop}

\bigskip

\begin{proof}
As the formula is additive, we may, and we will, assume that $E$ is an irreducible 
$\CM W$-module, concentrated in degree $0$. 

Now, let $\chi \in \Irr(W)$: we denote by $\mGB_\chi$ the maximal ideal of 
$\Zb_c$ corresponding to the fixed point $z_\chi$. We set $\pGB=\mGB_\chi \cap \Pb$: 
it does not depend on $\chi$ (it is the maximal ideal of $\Pb=\CM[V/W \times V^*/W]$ 
of functions which vanishes at $(0,0)$; see for instance~\cite[(14.2.2)]{BR}). 
By~(\ref{eq:local-chern}),  
$$\chern_c(E)=\sum_{\chi \in \Irr(W)} 
\dim^\gr((e\Hb_c \otimes_{\CM W} E)/\mGB_\chi(e\Hb_c \otimes_{\CM W} E)) e_\chi^W.$$
Now, let $\Zbov_c=\Zb_c/\pGB \Zb_c$ and $\Hbov_c=\Hb_c/\pG\Hb_c$. 
We set $\mGBov_\chi=\mGB_\chi/\pGB \Zb_c$. 
Then $\Hbov_c$ is a finite dimensional $\CM$-algebra (called the {\it restricted 
rational Cherednik algebra}) and again the bimodule $e\Hbov_c$ induces 
a Morita equivalence between $\Hbov_c$ and $\Zbov_c$. This implies that 
$e\Hbov_c/\mGBov_\chi \Hbov_c = (\Zbov_c/\mGBov_\chi) \otimes_{\Zbov_c} e\Hbov_c$ 
is a simple right $\Hbov_c$-module (which will be denoted by $\LC_c(\chi)$: it is isomorphic 
to the shift by some $r_\chi \in \ZM$ of the quotient of the baby Verma module denoted by 
$L(\chi)$ in~\cite{gordon}). 

Since $\CM W$ is semisimple, the $\CM W$-module $E$ is flat, and so 
\eqna
(e\Hb_c \otimes_{\CM W} E)/\mGB_\chi(e\Hb_c \otimes_{\CM W} E)
&\simeq& (e\Hb_c/\mGB_\chi(e\Hb_c)) \otimes_{\CM W} E \\
&\simeq& (e\Hbov_c/\mGBov_\chi(e\Hbov_c)) \otimes_{\CM W} E.\\
&\simeq& \LC_c(\chi) \otimes_{\CM W} E.
\endeqna
But the graded dimension of $\LC_c(\chi) \otimes_{\CM W} E$ is known 
whenever $\ZCB_{\! c}$ is smooth and is given by the expected formula 
(see~\cite[Lemma~3.3 and its proof]{bellamy}), up to a shift in grading:
$$\chern_c(E) =\sum_{\chi \in \Irr(W)} q^{r_\chi'}
\frac{\langle \chi , \CM[V]^\cow \otimes E \rangle_W^\gr}{\la \chi, \CM[V]^\cow \ra_W^\gr} e_\chi^W,$$
where $r_\chi' \in \ZM$ does not depend on $E$. Now, if $\CM$ denotes the trivial $\CM W$-module concentrated 
in degree $0$, then $\chern_c(\CM)=1$, which shows that $r_\chi'=0$ for all $\chi$, as desired.
\end{proof}

\bigskip

Let $\Krm_{\CM^\times}(\CM W)$ denote the Grothendieck group of the category of finite dimensional graded 
$\CM W$-modules. If $E$ is a finite dimensional graded $\CM W$-module, we denote 
by $\isomorphisme{E}$ its class in $\Krm_{\CM^\times}(\CM W)$. We still denote by 
$\chern_c : \Krm_{\CM^\times}(\CM W) \to \CM[[\hbar]] \otimes \Zrm(\CM W)$ 
the map defined by 
$$\chern_c(\isomorphisme{E})=\chern_c(E).$$
Now, let $W'$ be a parabolic subgroup of $W$ and set $V'=V^{W'}$ and $r=\codim_\CM(V')$. 
We identify the dual $V^{\prime *}$ of $V'$ with $V^{*W'}$ and note that 
\equat\label{eq:vprime}
V^*=V^{\prime *} \oplus (V^{\prime})^\perp.
\endequat
We denote by $\wedge (V^{\prime})^\perp$ the element of $\Krm_{\CM^\times}(\CM W')$ 
defined by
$$\wedge (V^{\prime})^\perp = \sum_{i \ge 0} (-1)^i \isomorphisme{\wedge^i (V^{\prime})^\perp}.$$
Recall also that there exist $n$ algebraically independent homogeneous polynomials 
$f_1$,\dots, $f_n$ in $\CM[V]^W$ such that $\CM[V]^W=\CM[f_1,\dots,f_n]$, 
and we denote by $d_i$ the degree of $f_i$. 

\bigskip

\begin{coro}\label{coro:chern}
Assume that $\ZCB_{\! c}$ is smooth. 
Let $E'$ be a finite dimensional graded $\CM W'$-module. Then
$$\chern_c\bigl(\Ind_{W'}^W(\wedge (V^{\prime})^\perp \otimes \isomorphisme{E'})\bigr) 
= \frac{(1-q^{d_1})\cdots(1-q^{d_n})}{(1-q)^{n-r}} 
\sum_{\chi \in \Irr(W)} 
\frac{ \langle \chi, \Ind_{W'}^W E' \rangle_W^\gr}{\la \chi, \CM[V]^\cow \ra_W^\gr } 
 e_\chi^W.$$
\end{coro}

\bigskip

\begin{proof}
The group $W'$ acts trivially on $V'$ so it acts trivially on $\wedge^i V^{\prime *}$ for all $i$. 
Therefore, 
$$(1-q)^{n-r}\chern_c\bigl(\Ind_{W'}^W(\wedge (V^{\prime})^\perp \otimes \isomorphisme{E'})\bigr)=
\chern_c(\Ind_{W'}^W\bigl(\wedge V^{\prime *} \otimes \wedge (V^{\prime})^\perp \otimes 
\isomorphisme{E'})\bigr).$$
But $\wedge V^{\prime *} \otimes \wedge (V^{\prime})^\perp=\Res_{W'}^W (\wedge V^*)$ by~(\ref{eq:vprime}), 
so, by Frobenius formula, 
$$(1-q)^{n-r}\chern_c\bigl(\Ind_{W'}^W(\wedge (V^{\prime})^\perp \otimes \isomorphisme{E'})\bigr) 
= \chern_c\bigl(\wedge V^* \otimes \Ind_{W'}^W E' \bigr).$$
So it follows from Proposition~\ref{prop:chern} that 
\begin{multline*}
(1-q)^{n-r}\chern_c\bigl(\Ind_{W'}^W(\wedge (V^{\prime})^\perp \otimes \isomorphisme{E'})\bigr) \\
= \sum_{\chi \in \Irr(W)} \frac{1}{\la \chi, \CM[V]^\cow \ra_W^\gr} 
\la \chi , \CM[V]^\cow \otimes (\wedge V^*) \otimes \Ind_{W'}^W(E') \ra_W^\gr ~e_\chi^W. 
\end{multline*}
But, if $w \in W$, the Molien's formula implies that the graded trace of $w$ on 
$\CM[V]^\cow$ is equal to 
$$\frac{(1-q^{d_1})\cdots(1-q^{d_n})}{\det(1-w^{-1} q)},$$
while its graded trace on $\wedge V^*$ is equal to $\det(1-w^{-1} q)$. So 
the class of $\CM[V]^\cow \otimes \wedge V^*$ is equal to $(1-q^{d_1})\cdots(1-q^{d_n})$  
times the class of the trivial module, and the corollary follows.
\end{proof}

\bigskip

\begin{proof}[Proof of Theorem~A]
Assume that $\ZCB_{\! c}$ is smooth. 
This implies in particular that Conjecture~\ref{conj:h}(1) and~(2) hold (Etingof-Ginzburg 
Theorem~\ref{theo:eg}), 
and so the hypotheses of Proposition~\ref{prop:reduction} are satisfied. Note also 
that $\im(\Omeb^c)=\Zrm(\CM W)$ by~(\ref{eq:lissite-omega}). It is then 
sufficient to prove that $\hbar^i \FC_i(\Zrm(\CM W)) \subset \im(\local_c^*)$  
for all $i$. 

Let us introduce some notation. If $G$ is a finite group and $H$ is a subgroup, 
we define a $\CM$-linear map $\Tr_H^G : \Zrm(\CM H) \longto \Zrm(\CM G)$ 
by 
$$\Tr_H^G(z)=\sum_{g \in [G/H]} \lexp{g}{z}=\frac{1}{|H|} \sum_{g \in G} \lexp{g}{z}$$
(here, $[G/H]$ denotes a set of representatives of elements of $G/H$ and $\lexp{g}{z}=gzg^{-1}$). 
It is easy to check that 
\equat\label{eq:trace-e}
\Tr_H^G(e_\eta^H) = \frac{\eta(1)}{|H|} \sum_{\g \in \Irr(G)} \frac{|G|}{\g(1)} 
\la \g,\Ind_H^G(\eta) \ra_G e_\g^G
\endequat
for all $\eta \in \Irr(H)$. Also, if $h \in H$ and $\Sigma_H(h) \in \Zrm(\CM H)$ denotes 
the sum of the conjugates of $h$ in $H$, then 
\equat\label{eq:trace-sigma}
\Tr_H^G(\Sigma_H(h))=\frac{|C_G(h)|}{|C_H(h)|} \Sigma_G(h),
\endequat
where $C_G(h)$ and $C_H(h)$ denote the centralizers of $h$ 
in $G$ and $H$ respectively.
Let $\PC_r(W)$ denote the set of parabolic subgroups $W'$ of $W$ such that 
$\codim_\CM(V^{W'}) = r$. It follows 
from~(\ref{eq:trace-sigma}) that 
\equat\label{eq:image-trace}
\FC_r(\Zrm(\CM W)) = \sum_{W' \in \PC_r(W)} \Tr_{W'}^W(\Zrm(\CM W')).
\endequat
Therefore, by Proposition~\ref{prop:reduction}, it is sufficient to prove that 
$$\text{\it $\hbar^r \Tr_{W'}^W(e_{\chi'}^{W'}) \in \im(\local_c^*)$ 
for all $W' \in \PC_r(W)$ and all $\chi' \in \Irr(W')$.}\leqno{(\bigstar)}$$
But it turns out that the coefficient of $\hbar^r$ in 
$$\frac{\chi'(1)}{|W'|}\chern_c\bigl(\Ind_{W'}^W(\wedge (V^{\prime})^\perp \otimes \chi')\bigr)$$
is equal, according to Corollary~\ref{coro:chern}, to
$$(-1)^r\frac{\chi'(1)}{|W'|}\sum_{\chi \in \Irr(W)} \frac{|W|~\la \chi, \Ind_{W'}^W \chi' \ra_W}{\chi(1)} ~e_\chi^W.$$
Indeed, since $q=\exp(\hbar)$, this follows from the fact that the polynomial $\la \chi, \CM[V]^\cow \ra_W^\gr$ 
takes the value $\chi(1)$ whenever $q=1$ (i.e. $\hbar =0$) and from the fact that
$$\frac{(1-q^{d_1})\cdots(1-q^{d_n})}{(1-q)^{n-r}} \equiv (-1)^r d_1\cdots d_n \hbar^r \mod \hbar^{r+1}$$
and that $|W|=d_1\cdots d_n$. Note that by definition 
$\frac{\chi'(1)}{|W'|}\chern_c\bigl(\Ind_{W'}^W(\wedge (V^{\prime})^\perp \otimes \chi')\bigr)$ belongs 
to $\CM[[\hbar]] \otimes_{\CM[\hbar]} \im(\local_c^*)$, so each of its homogeneous 
components belong to $\im(\local_c^*)$. 
By~(\ref{eq:trace-e}), this proves that 
$(\bigstar)$ holds, and the proof of Theorem~A is complete.
\end{proof}

\bigskip

\section{Proof of Theorem~B}\label{sec:proof-B}

\medskip

\springer{\noindent{\bf Hypothesis and notation.} 
{\it We assume in this section, and only in this section, that 
the symplectic quotient singularity $\ZCB_{\! 0}=(V \times V^*)/W$ admits 
a symplectic resolution $\XCB \to \ZCB_{\! 0}$.}}

\medskip

Recall~\cite[Theorem~1.3(ii)]{K} that the $\CM^\times$-action on $\ZCB_{\! 0}$ 
lifts uniquely to $\XCB$. 
As it is explained in Remark~\ref{rem:classification}, the existence of 
a symplectic resolution of $\ZCB_{\! 0}$ implies that 
all the irreducible components of $W$ are of type $G(d,1,n)$ or $G_4$. 
Since the proof of Theorem~B can be easily reduced to the irreducible case, 
we will separate the proof in two cases.

\begin{proof}[Proof of Theorem~B for $W=G(d,1,n)$]
Assume here that $W=G(d,1,n)$. 
Let $\Sb^1$ be the group of complex numbers of modulus $1$. 
In this case, it follows from~\cite{gordon o} that $\XCB$ 
is diffeomorphic to some smooth $\ZCB_{\! c}$, and that the diffeomorphism might be chosen to be 
$\Sb^1$-equivariant. As the $\Sb^1$-equivariant cohomology 
is canonically isomorphic to the $\CM^\times$-equivariant cohomology, 
this proves Theorem~B in this case.
\end{proof}

\begin{proof}[Proof of Theorem~B for $W=G_4$]
Assume here that $W=G_4$. It is possible (probable?) that again 
$\XCB$ is $\Sb^1$-diffeomorphic to some smooth $\ZCB_{\! c}$, 
but we are unable to prove it (it is only known that they 
are diffeomorphic). So we will prove Theorem~B in this case 
by brute force computations. 

We fix a primitive third root of unity $\z$ and we assume that 
$V=\CM^2$ and that $W=\langle s,t\rangle$, where 
$$s=\begin{pmatrix} \z & 0 \\ \z^2 & 1 \end{pmatrix} 
\qquad \text{and}\qquad 
t=\begin{pmatrix} 1 & -\z^2 \\ 0 & \z \end{pmatrix}.
$$
%We set for simplification $r=2\z+1$, so that $r^2=-3$. 
By the work of Bellamy~\cite{bellamy reso}, there are only two symplectic resolutions 
of $\ZCB_{\! 0}=(V \times V^*)/W$. They have both been constructed by Lehn and Sorger~\cite{LS}: 
one can be obtained from the other by exchanging the role of $V$ and $V^*$, so we will only 
prove Theorem~B for one of them. Let us describe it. 

Let $H=V^s$ and let $\HCB$ denote the image of $H \times V^*$ 
in $\ZCB_{\! 0}$, with its reduced structure of closed subvariety. 
We denote by $\b : \YCB \to \ZCB_{\! 0}$ the blowing-up of $\ZCB_{\! 0}$ along 
$\HCB$ and we denote by $\a : \XCB \to \YCB$ the blowing-up 
of $\YCB$ along its reduced singular locus $\SCB$. 
Then~\cite{LS}
$$\pi=\b \circ \a : \XCB \longto \ZCB_{\! 0}$$
is a symplectic resolution.

We now give more details, which all can be found in~\cite[\S{1}]{LS}. 
First, $\CM[\ZCB_{\! 0}]$ is generated 
by $8$ homogeneous elements $(z_i)_{1 \le i \le 8}$ whose degrees are given by the following table:
\equat\label{eq:deg-z}
\begin{array}{|c|c|c|c|c|c|c|c|c|}
\hline
z \vphantom{\DS{\frac{A}{A}}} & z_1 & z_2 & z_3 & z_4 & z_5 & z_6 & z_7 & z_8 \\
\hline
\deg(z) \vphantom{\DS{\frac{A}{A}}}& 0 & 4 & -4 & 2 & -2 & -6 & 6 & 0 \\
\hline
\end{array}
\endequat
The defining ideal of $\HCB$ in $\CM[\ZCB_{\! 0}]$ is generated by $6$ homogeneous elements 
$(b_j)_{1 \le j \le 6}$ whose degrees are given by the following table:
\equat\label{eq:deg-b}
\begin{array}{|c|c|c|c|c|c|c|}
\hline
b \vphantom{\DS{\frac{A}{A}}} & b_1 & b_2 & b_3 & b_4 & b_5 & b_6 \\
\hline
\deg(b) \vphantom{\DS{\frac{A}{A}}}& 2 & 6 & 0 & 12 & 8 & 4 \\
\hline
\end{array}
\endequat
This defines a $\CM^\times$-equivariant closed immersion $\YCB \injto \ZCB_{\! 0} \times \Pb^5(\CM)$. 
We denote by $\YCB_{\! i}$ the affine chart defined by ``$b_i \neq 0$''. The equations of the 
zero fiber $\b^{-1}(0)$ given in~\cite[\S{1}]{LS} show that $\YCB^{\CM^\times} = \{p_2,p_3,p_4,p_6\}$, 
where $p_i$ is the unique element of $\YCB_{\! i}^{\CM^\times}$. We use the notation of 
Example~\ref{ex:blowing}. By~(\ref{eq:local-blowing}) 
and~(\ref{eq:deg-b}), we have
$$\local_{\YCB}^*(\chern^1_\YCB(\isomorphisme{\b^*\HCB})) = 
-\hbar (6 \eb_{p_2} + 12 \eb_{p_4} + 4 \eb_{p_6}).\leqno{(\clubsuit)}$$
Now, $\SCB$ is contained in $\YCB_{\! 2} \cup \YCB_{\! 3}$, so $\a$ is an isomorphism 
in a neighborhood of $p_4$ and $p_6$. So let $q_4=\a^{-1}(p_4)$ and $q_6=\a^{-1}(p_6)$. 
These are elements of $\XCB^{\CM^\times}$. 

On the other hand, $\YCB_{\! 2}$ is a transversal $A_1$-singularity, so the defining 
ideal of $\SCB \cap \YCB_{\! 2}$ in $\CM[\YCB_{\! 2}]$ is generated by three homogeneous 
elements $a_+$, $a_\circ$ and $a_-$ (of degree $6$, $0$ and $-6$, 
by~\cite[\S{1}]{LS}), and it is easily checked 
that $\a^{-1}(p_2)^{\CM^\times}=\{q_2^+,q_2^-\}$ where $q_2^\pm$ is the unique $\CM^\times$-fixed 
element in the affine chart defined by ``$a_\pm \neq 0$''. 
\def\hilbert{{\mathbf{Hilb}}}

Also, $\YCB_{\! 3}$ is isomorphic to $(\hG \times \hG^*)/\SG_3$, where $\hG$ is the diagonal Cartan subalgebra 
of $\sG\lG_3(\CM)$ and $\SG_3$ is the symmetric group on $3$ letters, viewed 
as the Weyl group of $\sG\lG_3(\CM)$. Let $\hilbert_3(\CM^2)$ denote the Hilbert 
scheme of $3$ points in $\CM^2$, and let $\hilbert_3^0(\CM^2)$ denote the 
(reduced) closed subscheme defined as the Hilbert scheme of three points in $\CM^2$ 
whose sum is equal to $(0,0)$. By~\cite[Proposition~2.6]{haiman}, 
\equat\label{eq:x3}
\XCB_{\! 3} \simeq \hilbert_3^0(\CM^2).
\endequat
It just might be noticed that the isomorphism $\YCB_{\! 3} \simeq (\hG \times \hG^*)/\SG_3$ becomes 
$\CM^\times$-equivariant if one ``doubles the degrees'' in $(\hG \times \hG^*)/\SG_3$, that is, 
if $\CM^\times$ acts on $\hG$ (respectively $\hG^*$) with weight $2$ (respectively $-2$). 
Recall that $\CM^\times$-fixed points in $\hilbert_3^0(\CM^2)$ are parametrized 
by partitions of $3$: we denote by $q_3^+$, $q_3^\circ$ and $q_3^-$ the fixed points in 
$\XCB_{\! 3}$ corresponding respectively to the partitions $(3)$, $(2,1)$ and $(1,1,1)$ 
of $3$, so that
\equat\label{eq:x3-fixe}
\XCB_{\! 3}^{\CM^\times}=\{q_3^+,q_3^\circ,q_3^-\}.
\endequat
Finally, we have
\equat\label{eq:x-fixe}
\XCB^{\CM^\times}=\{q_2^+,q_2^-,q_3^+,q_3^\circ,q_3^-,q_4,q_6\}.
\endequat
Also, $\pi^*(\HCB)=\a^*(\b^*(\HCB))$ is an effective Cartier divisor of $\XCB$, 
and it follows from~($\clubsuit$) and the commutativity of the diagram~\eqref{eq:chern} that
$$\local_{\XCB}^*(\chern^1_\XCB(\isomorphisme{\pi^* \HCB})) = 
-\hbar (6 \eb_{q_2^+} + 6 \eb_{q_2^-} + 12 \eb_{q_4} + 4 \eb_{q_6}).\leqno{(\diamondsuit)}$$

We now wish to compute $\local_{\XCB}^*(\chern^1_\XCB(\isomorphisme{\a^* \SCB}))$. 
As the singular locus $\SCB$ is contained in $\YCB_{\! 2} \cup \YCB_{\! 3}$, and contains 
$p_2$ and $p_3$, there exists $n_2^+$, $n_2^-$, $n_3^+$, $n_3^\circ$ and $n_3^-$ in $\ZM$ such 
that
$$\local_{\XCB}^*(\chern^1_\XCB(\isomorphisme{\a^* \SCB}))=-\hbar(n_2^+ \eb_{q_2^+} + 
n_2^- \eb_{q_2^-} + n_3^+ \eb_{q_3^+} + n_3^\circ \eb_{q_3^\circ} + n_3^- \eb_{q_3^-}).$$
Since $a_+$ and $a_-$ have degree $6$ and $-6$, it follows from~(\ref{eq:local-blowing}) 
that $n_2^+=6$ and $n_2^-=-6$. 

As we can exchange the roles of $\hG$ and $\hG^*$ in the description of $\YCB_{\! 3}$ 
and its singular locus (because $\hG \simeq \hG^*$ as an $\SG_3$-module), this shows 
that $n_3^-=-n_3^+$ and $n_3^\circ=-n_3^\circ$. So $n_3^\circ=0$ and it remains 
to compute $n_3^+$. So let $U_+$ be the open subset of $\hilbert_3^0(\CM^2)$ 
consisting of ideals $J$ of codimension $3$ of $\CM[x,y]$ 
such that the classes of $1$, $x$ and $x^2$ form a basis 
of $\CM[x,y]/J$. Then we have an isomorphism $J_+ : \CM^4 \longiso U_+$ given by 
$$(a,b,c,d) \longmapsto J_+(a,b,c,d)=\langle x^3+ax+b,y- cx^2-dx-\frac{2}{3}ac \rangle.$$
The form of the generators of the ideal $J_+(a,b,c,d)$ is here to ensure that 
$J_+(a,b,c,d) \in \hilbert_3^0(\CM^2)$. 
The fixed point $q_3^+$ is the unique one in $U_+$. 
Through this identification, the action of $\CM^\times$ on $\CM^4$ is given 
by 
$$\xi\cdot(a,b,c,d)=(\xi^4 a, \xi^6 b , \xi^{-6} c ,\xi^{-4} d)$$
(remember that we must ``double the degrees'' of the usual action). 
The equation of $\a^*(\SCB)$ on this affine chart $\simeq \CM^4$ can then be computed 
explicitly and is given by $4 a^3 + 27 b^2=0$, so is of degree $12$. This shows that $n_3^+=12$. 
Finally,
$$\local_{\XCB}^*(\chern^1_\XCB(\isomorphisme{\a^* \SCB}))=-\hbar(6 \eb_{q_2^+} - 
6 \eb_{q_2^-} + 12 \eb_{q_3^+} -12 \eb_{q_3^-}).\leqno{(\heartsuit)}$$

Let us now conclude. First, recall from~\cite[Theorem~1.2]{GK} that 
\equat\label{eq:coho}
\begin{cases}
\Hrm^{2i+1}(\XCB) = 0 & \text{for all $i$,}\\
\Hrm^{2\bullet}(\XCB) \simeq \grad_\FC(\Zrm(\CM W)). 
\end{cases}
\endequat
By Proposition~\ref{prop:localisation}, we get 
$$\dim_\CM \Hrm^{2i}_{\CM^\times}(\XCB)=
\begin{cases}
1 & \text{if $i=0$,}\\
3 & \text{if $i=1$,}\\
7 & \text{if $i \ge 2$.}
\end{cases}\leqno{(\spadesuit)}$$
Now, $W$ has seven irreducible characters $1$, $\e$, $\e^2$, $\chi$, $\chi\e$, $\chi\e^2$ and $\th$, 
where $\chi$ is the unique irreducible character of degree $2$ with rational values and $\th$ 
is the unique one of degree $3$. We denote by 
$$\Psi : \Hrm^{2\bullet}_{\CM^\times}(\XCB^{\CM^\times}) \longiso \CM[\hbar] \otimes \Zrm(\CM W)$$
the isomorphism of $\CM[\hbar]$-algebras such that 
$$\Psi(\eb_{q_4})=e_1^W,\quad\Psi(\eb_{q_6})=e_\th^W,\quad\Psi(\eb_{q_2^+})=e_{\chi \e}^W,
\quad\Psi(\eb_{q_2^-})=e_{\chi\e^2}^W,$$
$$\Psi(\eb_{q_3^+}) = e_{\e^2}^W,\quad\Psi(\eb_{q_3^\circ})=e_\chi^W\quad\text{and}\quad
\Psi(\eb_{q_3^-})=e_\e^W.$$
By Proposition~\ref{prop:localisation}, $\Hrm^{2\bullet}_{\CM^\times}(\XCB)$ 
is isomorphic to its image by $\Psi \circ \local_{\XCB}^*$ in $\CM[\hbar] \otimes \Zrm(\CM W)$. 
But, by~($\diamondsuit$) and~($\heartsuit$), and after investigation of the character table of $W$, 
this image contains 
$$\hbar(6 e_{\chi\e}^W - 6 e_{\chi\e^2}^W + 12 e_1^W + 4 e_\th^W)=\hbar (4 + \Sigma_W(s) + \Sigma_W(s^2))$$
$$\hbar(12e_{\e^2}^W-12e_{\e}^W+6 e_{\chi \e}-6e_{\chi\e^2})=\hbar ((1+2\z) \Sigma_W(s) + (1+2\z^2)\Sigma_W(s^2)).
\leqno{\text{and}}$$
So it contains $\hbar \Sigma_W(s)$ and $\hbar \Sigma_W(s^2)$. Also, 
by~($\spadesuit$) and Proposition~\ref{prop:localisation}, it also contains 
$\hbar^2 \Zrm(\CM W)$, so 
$$\Rees_\FC^\bullet(\Zrm(\CM W)) \subset \im(\Psi \circ \local_{\XCB}^*).$$
Using again~($\spadesuit$) and Proposition~\ref{prop:localisation}, a comparison of 
dimensions yields that 
$$\Rees_\FC^\bullet(\Zrm(\CM W)) = \im(\Psi \circ \local_{\XCB}^*),$$
and the proof is complete.
\end{proof}

\end{document}